\documentclass{cocv}

\usepackage[english]{babel}
\usepackage[utf8]{inputenc}
\usepackage{amssymb}
\usepackage{graphicx}
\usepackage{color}
\usepackage{hyperref}
\hypersetup{pdftex,bookmarks=true}
\hypersetup{
pdfauthor = {Matthias Sonntag, Stephan Schmidt, Nicolas R. Gauger},
pdftitle = {Shape Derivatives for the Compressible Navier-Stokes Equations in Variational Form},
pdfsubject = {Shape optimization based on surface gradients and the Hadarmard-form is considered for a compressible viscous fluid. Special attention is given to the difference between the 'function composition' approach involving local shape derivatives and an alternate methodology based on the weak form of the state equation. The resulting gradient expressions are found to be equal only if the existence of a strong form solution is assumed. Surface shape derivatives based on both formulations are implemented within a Discontinuous Galerkin flow solver of variable order. The gradient expression stemming from the variational approach is found to give superior accuracy when compared to finite differences.},
pdfkeywords = {Shape derivative, Variational form, Compressible Navier-Stokes equations}}

\newtheoremstyle{ts}{}{}{\normalfont}{}{\normalfont\bfseries}{}{\newline}{\thmname{#1}\thmnumber{ #2}\thmnote{ (#3)}}
\theoremstyle{ts}
\newtheorem{thm}{Theorem}[section] 
\newtheorem{kor}[thm]{Corollary}
\newtheorem{de}[thm]{Definition}

\newtheorem{rem}[thm]{Remark}

\newenvironment{prf}[1][\proofname]{%
  \proof[\normalfont\bfseries #1:]%
}{\endproof}

\newcommand{\gradient}{\nabla}
\newcommand{\F}{\mathcal{F}}
\newcommand{\Fc}{\F^c(\u)}
\newcommand{\Fv}{\F^v(\u, \gradient \u)}
\newcommand{\Fcu}{\F^c_\u(\u)}
\newcommand{\Fvu}{\F^v_\u(\u, \gradient \u)}
\newcommand{\Fvgradu}{\F^v_{\gradient \u}(\u, \gradient \u)}
\newcommand{\fracp}[2]{\frac{\partial #1}{\partial #2}}
\renewcommand{\u}{\textbf{u}}
\newcommand{\z}{\textbf{z}}
\renewcommand{\v}{\textbf{v}}
\newcommand{\w}{\textbf{w}}
\newcommand{\dd}{\mbox{\,d}}
\newcommand{\definition}[1]{#1}
\renewcommand{\l}{\left}
\renewcommand{\r}{\right}
\newcommand{\transpose}{^{\top}}
\newcommand{\inside}{\text{ in }}
\newcommand{\on}{\text{ on }}
\newcommand{\Rd}{\mathbb{R}^d} 
\newcommand{\Dk}{C^k_0(D; \Rd)}
\newcommand{\dJ}{dJ(\Omega; V)}
\newcommand{\Vn}{\langle V, n \rangle}
\renewcommand{\div}{\operatorname{div}}
\newcommand{\ddt}[1]{\l. \frac{d #1}{dt} \r|_{t=0}}
\newcommand{\dV}[1]{d_V \! \left[#1\right]}
\newcommand{\sg}[2]{\l(#1, #2\r)_\Gamma}  
\newcommand{\sgw}[2]{\l(#1, #2\r)_{\Gamma_W}}
\newcommand{\sgf}[2]{\l(#1, #2\r)_{\Gamma \setminus \Gamma_W}}
\newcommand{\so}[2]{\l(#1, #2\r)_\Omega}    
\newcommand{\s}[3]{\l(#1, #2\r)_{#3}}
\newcommand{\extrabox}[1]{\boxed{#1}} 
\newcommand{\num}[1]{\text{\fboxsep0.2em\fbox{$#1$}}}

\def\bei#1{\vrule width 0.4pt height 14pt depth 9pt
           \lower 8pt \hbox{$ _{\hbox{}\, #1}$}\!}

\begin{document}
\title{Shape Derivatives for the Compressible Navier-Stokes Equations in Variational Form}
\date{20.12.2013}
\author{Matthias Sonntag}\address{RWTH Aachen University, Aachen, Germany; \email{Sonntag@mathcces.rwth-aachen.de \& gauger@mathcces.rwth-aachen.de}}
\author{Stephan Schmidt}\address{Imperial College London, London, United Kingdom; \email{s.schmidt@imperial.ac.uk}}
\author{Nicolas R. Gauger}\sameaddress{1}

\begin{abstract}
Shape optimization based on surface gradients and the Hadarmard-form is considered for a compressible viscous fluid. Special attention is given to the difference between the ``function composition'' approach involving local shape derivatives and an alternate methodology based on the weak form of the state equation. The resulting gradient expressions are found to be equal only if the existence of a strong form solution is assumed. Surface shape derivatives based on both formulations are implemented within a Discontinuous Galerkin flow solver of variable order. The gradient expression stemming from the variational approach is found to give superior accuracy when compared to finite differences.
\end{abstract}
\maketitle
\section*{Introduction}
Shape optimization is a research field that has received much attention in the past. In general, any problem where the boundary of the domain is part of the unknown can be considered a shape optimization problem. In most applications, the physics are modelled by partial differential equations, making shape optimization a special sub-class within the field of PDE-constrained optimization. Usually, the derivation of the sensitivities and adjoint equations follows a function composition approach, i.e.\ some set of design variables defines the geometry and within this geometry the PDE is solved, thereby generating the state variables that enter the objective function~\cite{giles1997adjoint, Jameson88, PironneauxX2001}. Therefore, the necessity to consider sensitivities or derivative information with respect to the geometry adds additional complexity to the shape optimization problem when compared to general PDE-constrained optimization. Because it is often not immediately clear how to compute these ``mesh sensitivities'', that is the variation of the PDE with respect to a change in the geometry, there is often a strong desire for a very smooth parametrization of the domain with as few design parameters as possible. Although there have been successful attempts to incorporate problem structure exploitations in order to efficiently compute these partial derivatives for very large problems, such as differentiating the entire design chain at once or by considering the adjoint process of the mesh deformation~\cite{GaWaMoWi07, NielsenPark2006}, very often one is still forced into finite differencing, which means the PDE residual at steady state has to be evaluated on meshes that have been perturbed by a variation in each design parameter of the shape, a process that makes large scale optimization usually prohibitive. This negates some of the advantages of the adjoint approach, such as the independence of the number of design parameters. More severely, it also makes fast optimization strategies such as the one-shot approach~\cite{GriewankxX2005, GS09, TaAsanSxX1992a} somewhat unattractive in terms of wall-clock-time.

A more recent trend to overcome the cumbersome computation of these geometric sensitivities is the use of shape calculus. Shape calculus summarizes the mathematical framework used when considering problems where the shape is the unknown in the continuous setting. Manipulations in the tangent space of the unknown object can be used to circumvent any necessity of knowing discrete geometric sensitivities, because these can be directly included in a surface gradient expression on the continuous level. More details on this theoretical framework can be found in~\cite{shapesandgeometries, sokolowski}. Traditionally, this methodology was primarily used to address the very difficult question of existence and uniqueness of optimal shapes~\cite{Piro73}, but more recently this methodology has also been used in very large scale aerodynamic design and computational optimization~\cite{castro, SIGS2011, SS10GeneralNS}. In~\cite{SSIG2010A}, for example, the complete optimization of a blended wing-body aircraft in a compressible inviscid fluid is considered. Because this approach solely relies on the problem formulation in the continuous setting and only afterwards discretizes the continuous boundary integral expressions for the shape derivative, great care must be taken when making the initial assumptions and when implementing the respective continuous expressions, especially at singular points in the geometry, such as the trailing edge of an airfoil~\cite{Lozano2012}. Because this approach is indeed truly independent of the number of design parameters, it enables the most detailed possible parameterization, that is using all surface mesh nodes as a design unknowns. This is sometimes called ``free node parameterization''. However, these highly detailed shape parameterizations usually lack any kind of inherent regularity preservation and as such, one usually finds this approach paired with some sort of smoothing procedure that projects or embeds the respective optimization iterations into a desired regularity class, which can nicely be paired with an SQP or Newton-type optimization scheme, which is sometimes also called a ``Sobolev Method''~\cite{Arian2, schmidt:2562}.

As part of this work, we study how to further increase the accuracy of shape derivatives when used within viscous compressible aerodynamic design optimization. Within applied aerodynamic shape optimization, it is customary to exploit the above mentioned function composition approach in order to derive and implement the adjoint equation and gradient expression. This has been used with great success, both within the context of continuous and discrete adjoint based aerodynamic shape optimization~\cite{Zingg:AnaheimJournal, Papdimitriou2007, SIGS2011, Giannakoglou2009B} and general shape optimization~\cite{Eppler2007A}. However, common to these approaches is the assumption that the state equation has a strong form solution and each of the steps within the shape differentiation process of the function composition exists, which usually means the existence of so-called local shape derivatives. For elliptic problems, this existence can usually be shown, making the above mentioned approach somewhat of an established procedure, see for example Chapter 3.3 in~\cite{sokolowski}. However, for the hyperbolic equations governing some compressible fluids, the existence of a strong form solution is not clear. Rather, in the presence of shock waves and discontinuities in the flow, one can usually only expect the variational form of the equation to hold, a property which is very often not taken into account when studying the derivative. Shape differentiation of problems governed by PDEs in weak or variational forms are not very often considered in the literature, except in~\cite{ItoKunischPeichel2006} and especially in~\cite{KunishVariationalSD}, where the incompressible Navier--Stokes equations are considered for this purpose from a rigorous theoretical standpoint. Thus, we revisit the shape optimization problem previously considered in~\cite{Soemarwoto1}, but the gradient is derived using elements of the variational approach as shown in~\cite{KunishVariationalSD}. Furthermore, we simultaneously follow the function composition approach, outlining the exact differences comparing these two approaches. One can nicely see how both methodologies reduce to the same gradient expression when assuming the existence of a strong from solution of the state equation. We conclude with a numerical error analysis based on comparing finite differencing with either implementation, demonstrating the higher accuracy of the gradient formulation based on the variational form of the compressible Navier--Stokes equations.

The structure of the paper is as follows. In Section~\ref{sec:FluidMechanics}, we begin by recapitulating the compressible Navier--Stokes equations in both strong and variational form. Next, Section~\ref{section:shape_calculus} serves as an introduction and quick overview of shape calculus, including shape derivatives and the Hadamard or Hadamard--Zol\'esio Structure Theorem, which leads to a preliminary form of the shape derivative of the aerodynamic cost functions. The next section, Section~\ref{sec:NSSD}, is used to work out the differences between the shape derivative of the compressible Navier-Stokes equations stemming from either the function composition or the variational approach. In Section~\ref{section:adjoint_calculus}, we then summarize the idea of adjoint calculus. This is used to differentiate the Navier--Stokes equations, thereby discussing the Hadamard form of the respective objective functions both for the strong as well as the variational form of the state constraint. Finally, in the last section, numerical results achieved with both methods are compared to shape derivatives computed by finite differences, showing a considerable gain in accuracy when using shape derivatives based on the variational form. 

\section{Fluid Mechanics}\label{sec:FluidMechanics}

\subsection{Flow domain and boundary conditions}
In the following $\rho$, $v = (v_1, v_2)\transpose$, $p$, $E$ and $T$ denote the density, velocity, pressure, total energy and temperature. 
The domain of the fluid is denoted by $\Omega$, with wall and far-field boundaries $\Gamma_W$ and $\Gamma_\infty$.
At the wall $\Gamma_W$, the no-slip boundary condition $v = 0$ is imposed for the velocity. With respect to temperature, either the isothermal boundary condition $T = T_W$ or the adiabatic boundary condition $\gradient T \cdot n = 0$ holds. The isothermal and adiabatic parts of the wall are named $\Gamma_{iso}$ and $\Gamma_{adia}$ and we assume $\Gamma_W = \Gamma_{iso} \cup \Gamma_{adia}$ disjoint.

Furthermore $\kappa$, $e$, $H$, $\mu$ and $\gamma$ denote the thermal conductivity, the internal energy, the enthalpy, the viscosity and the adiabatic exponent. The relation $T \kappa = \frac{\mu \gamma}{Pr} \l(E - \frac{1}{2} \Vert v\Vert^2 \r)$ is fulfilled for the temperature, where $Pr$ is the Prandtl number.

\subsection{Navier--Stokes equations and aerodynamic objective functions}
In this subsection, we state the Navier--Stokes equations in both strong and weak form. As discussed later, the shape derivative of the aerodynamic cost functions differs depending on which form of the Navier--Stokes equations is used.
Using the viscous stress tensor $\tau$, defined by 
\begin{equation}
\label{eq:ns_integral}
\tau = \mu \l( \gradient v + (\gradient v) \transpose - \frac{2}{3} (\gradient \cdot v) I \r)
,
\end{equation}
the compressible Navier-Stokes equations in strong form are given by
\begin{equation}
\label{eq:ns_pointwise}
\gradient \cdot \l( \Fc - \Fv \r) \equiv \sum_k \l( \fracp{}{x_k} f_k^c(\u) - \fracp{}{x_k}f_k^v(\u, \gradient \u) \r) = 0 \quad \inside \Omega,
\end{equation}
where $\u$ denotes the vector of conserved variables, $\F^c = (f^c_1, f^c_2)$ the convective fluxes
\begin{equation*}
\u = \begin{pmatrix}
\rho \\
\rho v_1 \\
\rho v_2 \\
\rho E
\end{pmatrix}
\quad \quad
f_1^c(\u) = \begin{pmatrix}
\rho v_1 \\
\rho v_1^2 + p \\
\rho v_1 v_2 \\
\rho H v_1
\end{pmatrix},
 \quad \quad 
f_2^c(\u) = \begin{pmatrix}
\rho v_2 \\
\rho v_1 v_2 \\
\rho v_2^2 + p \\
\rho H v_2
\end{pmatrix}
\end{equation*}
and $\F^v = (f^v_1, f^v_2)$ denotes the viscous fluxes
\begin{equation*}
f_1^v(\u, \gradient \u) = \begin{pmatrix}
0 \\
\tau_{11} \\
\tau_{21} \\
\sum_j \tau_{1j} v_j + \kappa \fracp{T}{x_1}
\end{pmatrix},
 \quad \quad 
f_2^v(\u, \gradient \u) = \begin{pmatrix}
0 \\
\tau_{12} \\
\tau_{22} \\
\sum_j \tau_{2j} v_j + \kappa \fracp{T}{x_2}
\end{pmatrix}
.
\end{equation*}
Furthermore, temperature $T$ and pressure $p$ are linked to the state variables using the perfect gas assumption, that is
\begin{align*}
 p = \rho RT = \rho R \frac{E - \frac{1}{2}\Vert v\Vert^2}{c_v} = \frac{R}{c_v}\rho\left(E - \frac{1}{2}\Vert v\Vert^2\right),
\end{align*}
where $e = c_v T$, $E = e + \frac{1}{2}\Vert v\Vert^2$ and $c_v$ denotes the heat capacity of the gas at constant volume.

Furthermore, the compressible Navier-Stokes equations in variational form are given by the following

\begin{de}[Variational Form of the Navier-Stokes equations]\label{Def:Weak}
We assume that $\mathcal{H} := H^1 \times H^2 \times H^3 \times H^4$, where $H^i$ is a suitable Hilbert-Space that is provided by the respective user given flow discretization. Multiplication of the pointwise Navier-Stokes~\eqref{eq:ns_pointwise} with an arbitrary test function $\v\in\mathcal{H}$ and integration by parts results in the problem to find $\u \in \mathcal{H}$, such that
\begin{equation}
\label{eq:ns_variational}
  \langle F(\u, \Omega), \v\rangle_{\mathcal{H}^*\times \mathcal{H}} := 
 -\so{\Fc - \Fv}{\gradient \v} + \sg{n \cdot \l(\Fc - \Fv\r)}{\v}  = 0 \quad \forall \v \in \mathcal{H} 
\end{equation}
with the boundary conditions 
\begin{equation}
\label{eq:bc}
\begin{aligned}
\begin{pmatrix} v_1 \\ v_2 \end{pmatrix} &= 0 &&\on \Gamma_W, \\
 \gradient T \cdot n &= 0 &&\on \Gamma_{adia}, \\
 T &= T_{wall} &&\on \Gamma_{iso}.
\end{aligned}
\end{equation}
\end{de}

\begin{de}[Cost Function]
The cost functions under consideration are the lift and drag coefficients, given by
\begin{equation}\label{eq:CostFunction}
J(\u) = \frac{1}{C_\infty} \int_{\Gamma_W} \l( p n - \tau n \r) \cdot \psi \dd s
,
\end{equation}
where $C_\infty$ is a constant and $\psi$ is either $\psi_l = (-\sin(\alpha), \cos(\alpha)) \transpose$ for the lift or $\psi_d = (\cos(\alpha), \sin(\alpha)) \transpose$ for the drag coefficient.
\end{de}

\section{Shape Calculus}
\label{section:shape_calculus}
\subsection{Definition of the shape derivative and the Hadamard Theorem}

In this section the concept of shape derivatives and especially the Hadamard Theorem, as stated in \cite{shapesandgeometries, sokolowski}, are summarized. Let therefore $D$, the so-called hold-all, be an open set in $\Rd$ and the domain $\Omega$ be a measurable subset of $D$. For vector fields $V \in \Dk$ the pertubation of identity 
\begin{equation*}
T_t[V] : D\times [0, \delta) \to \Rd, \quad (x, t) \mapsto x + t V(x)
\end{equation*}
is a common approach to describe a deformation $\Omega_t = T_t[V](\Omega) $ of the domain $\Omega$. With such a deformation, the shape derivative of a domain functional $J(\Omega)$ at $\Omega$ in the direction of a vector field $V \in \Dk$ is defined as the Eulerian derivative 
\begin{equation*}
\dJ := \lim_{t \searrow 0} \frac{J(\Omega_t) - J(\Omega)}{t}
.
\end{equation*}
The functional $J(\Omega)$ is called \definition{shape differentiable} at $\Omega$ if this Eulerian derivative exists for all directions $V$ and the mapping $G(\Omega) : \Dk \to \mathbb{R}, \ V \mapsto \dJ$ is linear and continuous.

If the functional $J(\Omega)$ is shape differentiable on measurable sets $\Omega \subset D$, then there exists the shape gradient $G(\Omega) \in \l(\Dk\r)^*$ such that
\begin{equation}
\label{eq:shape_calculus1}
\dJ = \langle G(\Omega), V \rangle_{\l(\Dk\r)^* \times \Dk} \quad \forall \ V \in \Dk
.
\end{equation}
This means that the shape derivative, as a directional derivative in direction $V$, can be computed by the dual pair, which is a generalized scalar product, of the shape gradient $G(\Omega)$ and the vector field $V$. 
Furthermore in \cite{sokolowski} it is shown, that if a vector field $V$ fulfills $V \cdot n = 0$, meaning it is tangential to the boundary $\Gamma = \partial \Omega \in C^k$, then the shape derivative in this direction vanishes, which is intuitively clear, since a deformation of the domain $\Omega$ in a tangential direction, being a form of reparameterization, does not change the domain and therefore the shape derivative vanishes, because the flow solution and the cost functional stay the same. 

In addition, there exists a continuous linear mapping $dJ(\Gamma; \cdot) : C^k(\Gamma) \to \mathbb{R}$ such that for all vector fields $V \in C^k(\overline D; \Rd)$ the relation
\begin{equation*}
\dJ = dJ(\Gamma; V \cdot n)
\end{equation*}
holds, which implies that the shape derivative only depends on the normal component of the vector field $V$ at the boundary of the domain. The following theorem states the existence of a scalar distribution $g(\Gamma)$, which takes the role of the shape gradient $G(\Omega)$ in equation \eqref{eq:shape_calculus1} for the mapping $dJ(\Gamma; \cdot)$ on the boundary of the domain.

\begin{thm}[Hadamard Theorem, Hadamard formula]
\label{thm:hadamard}
For every domain  $\Omega \subset D$ of class $C^k$, let $J(\Omega)$ be a shape differentiable function. Furthermore let the boundary $\Gamma$ be of class $C^{k-1}$. There exists the following scalar distribution $g(\Gamma) \in C^k_0(\Gamma)^*$, such that the shape gradient $G(\Omega) \in C^k_0(\Omega, \Rd)^*$ of $J(\Omega)$ is given by
\begin{equation*}
G(\Omega) = \gamma_\Gamma^*(g \cdot n),
\end{equation*}
where $\gamma_\Gamma \in L\l(C^k_0(\Omega, \Rd), C^k_0(\Gamma, \Rd)\r)$ and $\gamma_\Gamma^*$ denote the trace operator and its adjoint operator. In this situation, one can show that~\cite{sokolowski}
\begin{equation*}
\dJ = dJ(\Gamma; V \cdot n) = \langle g, V \cdot n \rangle_{\l(C^k_0(\Gamma, \Rd)\r)^* \times C^k_0(\Gamma, \Rd)}
.
\end{equation*}
If $g(\Gamma)$ is integrable over $\Gamma$, than the \definition{Hadamard Formula} 
\begin{equation*}
\dJ = \int_\Gamma (V\cdot n) g \dd s
\end{equation*}
is fulfilled. In the following we will call terms that are of the structure ``$(V \cdot n) \ldots$'' to be in \definition{Hadamard form}.
\end{thm}
\begin{prf}
A proof can for example be found in~\cite{shapesandgeometries} or in~\cite{sokolowski}.
\end{prf}

\begin{rem}
\label{rem:normal_direction}
Assuming the boundary $\Gamma_W$ is of sufficient regularity such that the tangential component 
only describes a reparameterization but no actual change of the shape, then the shape derivative $\dJ$ only depends on the normal component $(V \cdot n)n$ of the vector field $V$, which will later eliminate certain normal components within the derivation of the surface gradient expression of the cost functions. Further studies of the ramifications of this assumption can be found in~\cite{Lozano2012}.
\end{rem}

\begin{de}[Material derivative / Local shape derivative]
\label{de:material_derivative}
\label{de:local_shape_derivative}
The total derivative
\begin{equation*}
\dV{f}(x) := \ddt{} f(t, T_t(x))
\end{equation*}
of $f$ is called the \definition{material derivative}. Furthermore, the partial derivative 
\begin{equation*}
f'(x) := f'[V](x) := \frac{\partial}{\partial t} f(t, x)
\end{equation*}
is called the \definition{local shape derivative} of $f$.
\end{de}

\begin{rem}
\label{rem:material_derivative}
The material and the local shape derivative are related to each other by the chain rule if both exist, i.e.\
\begin{equation*}
\dV{f}(x) = f'[V](x) + \gradient f(0, x) \cdot V(x) = f' + \gradient f \cdot V
,
\end{equation*}
where we used $\ddt{} T_t(x) = \ddt{} \l( x + t V(x)\r) = V(x)$ for the perturbation of identity. If the geometry is as such that the local shape derivative does not exist, a sharp convex corner for example, one usually accepts the above formula as a definition of the local shape derivative instead.
\end{rem}

\subsection{Tangential calculus}
Before we state the shape derivative of volume and boundary integrals, we give a minimal summary of tangential calculus, which will later be used to derive a preliminary shape derivative of the lift and drag coefficients. This is also required to find the local shape derivative of quantities fulfilling Neumann boundary conditions. A more detailed discussion on tangential calculus can be found in~\cite{shapesandgeometries}, Chapter 8, Section 5.

Let $f\in C^1(\Gamma)$ be a function with a $C^1$-extension $F$ into a tubular neighbourhood of $\Gamma$. The \definition{tangential gradient} is the ordinary gradient minus its normal component, i.e.\ 
\begin{align}
\gradient_\Gamma f := \gradient F|_\Gamma - \fracp{F}{n} n.
\end{align}
Analogously, the \definition{tangential divergence} of a smooth vector field $W \in (C^1(\Gamma))^d \cap (C^1(\Omega))^d$ is defined by
\begin{align}\label{eq:TangentDiv}
\div_\Gamma W := \div W - DW n \cdot n.
\end{align}
For $f$ and $W$ as defined above, the \definition{tangential Green's formula} is given by
\begin{equation}
\label{eq:greens_formula}
\int_\Gamma W \cdot \gradient_\Gamma f \dd s = \int_\Gamma f K (W \cdot n) - f \div_\Gamma W \dd s
,
\end{equation}
where $K:=\div_\Gamma n$ denotes the sum of the principal curvatures, the so called additive curvature, or $(d-1)$ times the mean curvature.
A proof can be found in \cite{shapesandgeometries}, Chapter 8. 

\subsection{Shape derivative for volume and boundary integrals}
In this subsection, we recapitulate shape derivatives of general volume and boundary integrals, as they can be found in \cite{shapesandgeometries} for example. Furthermore, we show a preliminary shape derivative of the drag and lift coefficients, which is later transformed into Hadamard form using adjoint calculus.
The shape derivative of a general volume cost function $J(\Omega) = \int_\Omega f(x) \dd x$ is given by
\begin{equation}
\label{eq:sd_volume}
\dJ = \int_\Omega f' \dd x + \int_\Gamma (V \cdot n) f \dd s
\end{equation}
and for a general boundary cost function $J(\Omega) = \int_\Gamma f \dd s$ the shape derivative fulfills 
\begin{equation}
\label{eq:sd_boundary}
\dJ = \int_\Gamma f' + (V \cdot n) \l(\fracp{f}{n} + K f \r)  \dd s
.
\end{equation}

Furthermore, if the vector field $V$ is orthogonal to the boundary $\Gamma$, the material derivative of the normal vector fulfils
\begin{equation}
\label{eq:sd_normal}
\dV{n} = - \gradient_\Gamma( V \cdot n)
, 
\end{equation}
which can be found in \cite{schmidt_thesis}.

\begin{thm}[Preliminary shape derivative of the cost functional]
\label{thm:sd_aerodynamic_targetfunction}
\label{thm:preliminary_sd}
If the vector field of the pertubation of identity fulfills $V = 0$ in the neighbourhood of the farfield boundary $\Gamma_\infty$, then the shape derivative of the lift and drag coefficients, Equation~\eqref{eq:CostFunction}, is given by
\begin{equation}\label{eq:PrelimGradient}
\dJ = \frac{1}{C_\infty} \int_{\Gamma_W} (p' n - \tau' n ) \cdot \psi  + (V \cdot n) \div(p\psi - \tau \psi) \dd s
.
\end{equation}
\end{thm}
\begin{prf}
The proof of the above expression is not fully straight forward, because the objective function~\eqref{eq:CostFunction} depends on the local normal, which first needs to be extended into a tubular neighbourhood as discussed within the subsection on tangential calculus. Following the same argumentation as in~\cite{schmidt_thesis}, let $\mathcal{N}$ be an extension of the normal $n$ into a tubular neighbourhood. It is easy to see that as a result of the normalization of this extension, the property
\begin{align}
 0 &= \nabla \mathcal{N} \cdot \mathcal{N} \nonumber \\
  &= \nabla \mathcal{N} \cdot n \text{ on } \Gamma_{W} \label{eq:normalex}
\end{align}
holds. Because $p\mathcal{N} \cdot \psi$ and $\tau \mathcal{N} \cdot \psi$ have the same structure, we restrict ourselves to $p \mathcal{N} \cdot \psi$, that is we consider the functional $j = \int_{\Gamma_W} p \mathcal{N} \cdot \psi \dd s$, for which the shape derivative of a general boundary integral as stated in equation \eqref{eq:sd_boundary} is applicable. Paired with the properties of the normal extension~\eqref{eq:normalex}, this leads to
\begin{align}
dj(\Omega; V) &= \int_{\Gamma_W} p' \mathcal{N} \cdot \psi + p \mathcal{N}' \cdot \psi + (V \cdot \mathcal{N}) \l( \fracp{(p\mathcal{N} \cdot \psi)}{n} + K (p\mathcal{N} \cdot \psi)  \r) \dd s \nonumber \\
&= \int_{\Gamma_W} p'n \cdot \psi + p n' \cdot \psi + (V \cdot n) \l( \fracp{(p\cdot \psi)}{n}n + K (pn \cdot \psi)  \r) \dd s \label{eq:proof_pre_sd_targetfunction}
.
\end{align}
From Remark \ref{rem:material_derivative} we have $\mathcal{N}' = \dV{\mathcal{N}} - \gradient \mathcal{N} \cdot V$ for the local shape derivative of the extended normal vector. Combining Remark~\ref{rem:normal_direction}, i.e.\ assuming $V=(V\cdot n)n$, with Equation~\eqref{eq:normalex}, one arrives at $n' = \dV{n}$. Also using $\dV{n} = - \gradient_\Gamma( V \cdot n)$ from \eqref{eq:sd_normal} gives us the relation
\begin{equation*}
\int_{\Gamma_W} p n' \cdot \psi \dd s = - \int_{\Gamma_W} p \gradient_\Gamma(V \cdot n) \cdot \psi \dd s
.
\end{equation*}
Application of the tangential Green's formula to the right hand side of the above equation yields
\begin{equation*}
\int_{\Gamma_W} p n' \cdot \psi \dd s = -\int_{\Gamma_W} (V \cdot n) K (p \psi \cdot n) - (V \cdot n) \div_\Gamma (p \psi) \dd s,
\end{equation*}
where the assumption $V = 0$ in the neighbourhood of the farfield boundary $\Gamma_\infty$ was used to employ the tangential Green's formula at the wall boundary only.
We will now insert the above equation into \eqref{eq:proof_pre_sd_targetfunction}. As one can see, the terms containing the additive curvature $K$ cancel out and the following expression remains
\begin{equation*}
dj(\Omega; V) = \int_{\Gamma_W} p' n \cdot \psi + (V \cdot n) \l[ \div_\Gamma(p \psi) + \fracp{(p \cdot \psi)}{n}n  \r] \dd s.
\end{equation*}
The terms within the bracket now exactly align with the definition of the tangential divergence, Equation~\eqref{eq:TangentDiv}, such that 
\begin{align*}
\div_\Gamma(p \psi) + \fracp{(p \cdot \psi)}{n}n &= \div_\Gamma(p \psi) + D(p \psi) n \cdot n = \div (p \psi).
\end{align*}
The same argumentation can now be also applied to $\tau n \cdot \psi$ instead of $p n \cdot \psi$.
\end{prf}

The above theorem already supplies one possible representation of the shape derivative of the objective functional. However, this preliminary shape derivative is not yet in Hadamard form because it still contains the local shape derivatives $p'$ and $\tau'$. Computation of these would require one forward flow solution for each design parameter of the parameterization of the shape, which is prohibitively costly. In Section \ref{section:adjoint_calculus}, adjoint calculus is used to remove these local shape derivatives $p'$ and $\tau'$, thereby transforming the above gradient expression into the Hadamard form.

\subsection{Shape derivatives of boundary conditions}
\label{section:sd_of_boundary_conditions}
Transformation of the gradient expression from Theorem~\ref{thm:sd_aerodynamic_targetfunction} requires explicit knowledge of the boundary conditions defining the local shape derivatives $p'$ and $\tau'$. These are governed by the respective boundary conditions imposed within the forward problem. As such, we now consider how the Dirichlet condition, the Neumann condition and the slip condition of the forward problem determine these local shape derivatives. The general argumentation again follows~\cite{sokolowski}.

\subsubsection{Dirichlet boundaries}\label{subsec:DirichletBC}
First, we consider a general Dirichlet boundary condition $w = w_D$ on the wall $\Gamma_W$, where $w_D$ does not depend on the geometry, meaning $w_D$ is independent of the parameter $t$ of the deformation $T_t$. This is especially true for the no-slip condition $v=0$, which should also be fulfilled on the perturbed boundary. Application of the material derivative paired with remark \ref{rem:material_derivative} applied to both sides of the Dirichlet boundary condition $w = w_D$ yields
\begin{equation*}
\dV{w} = w' + \gradient w \cdot V = \dV{w_D} = \gradient w_D \cdot V
,
\end{equation*}
where the local shape derivative of $w_D$ equals zero, since $w_D$ does not depend on $t$. From this equation, we can extract a condition defining the local shape derivative $w' = \gradient (w_D - w) \cdot V$.
According to Remark \ref{rem:normal_direction}, it is sufficient to consider only the normal direction $(V \cdot n) n$ of the vector field $V$, which leads to
\begin{equation*}
w' = \gradient (w_D - w) \cdot n(V \cdot n) = \fracp{w_D - w}{n} (V \cdot n).
\end{equation*}

\subsubsection{Neumann boundaries}
Similar to the Dirichlet boundary condition, we again would like to consider the material derivative of the boundary condition of the forward problem and then apply the chain rule argument given by Remark~\ref{rem:material_derivative} to find a corresponding expression for the local shape derivatives. This will again require the quantities under consideration to exist at least within a tubular neighbourhood for which we again assume an extension of the normal $\mathcal{N}$ just as in Theorem~\ref{thm:sd_aerodynamic_targetfunction}. If we apply the material derivative to both sides of the Neumann boundary condition $\fracp{w}{\mathcal{N}} = \gradient w \cdot \mathcal{N} = w_N$, where $w_N$ does not depend on the shape, meaning on the parameter $t$, we get with Remark \ref{rem:material_derivative} 
\begin{equation*}
\dV{\gradient w \cdot \mathcal{N} } = (\gradient w \cdot \mathcal{N})' + \gradient ( \gradient w \cdot \mathcal{N}) \cdot V = \dV{w_N} = \gradient w_N \cdot V
.
\end{equation*}
Since  $\gradient w' \cdot \mathcal{N} = (\gradient w \cdot \mathcal{N})' - \gradient w \cdot \mathcal{N}'$ holds, we get from this equation
\begin{equation*}
\gradient w' \cdot \mathcal{N} = \gradient w_N \cdot V - \gradient w \cdot \mathcal{N}' - \gradient ( \gradient w \cdot \mathcal{N}) \cdot V
.
\end{equation*}
Using the usual orthogonality argumentation again, Remark~\ref{rem:normal_direction} and again employing Remark~\ref{rem:material_derivative}, we can insert the relation $n' = \mathcal{N}' = \dV{\mathcal{N}} - \gradient \mathcal{N} \cdot V = \dV{n}$ and furthermore use  $\gradient(\gradient w \cdot \mathcal{N}) \cdot V = D^2 w \mathcal{N} \cdot V + \gradient w \cdot (\gradient \mathcal{N} \cdot V) = D^2 w n \cdot V$ to obtain
\begin{equation*}
\begin{aligned}
\gradient w' \cdot n &= \gradient w_N \cdot V - \gradient w \cdot \dV{n} - D^2 w n \cdot V \\
  &= (V \cdot n) \l[  \gradient w_N \cdot n - D^2w n \cdot n \r] - \gradient w \cdot \dV{n}.
\end{aligned}
\end{equation*}
From Equation~\eqref{eq:sd_normal} we get the equality $\dV{n} = -\gradient_\Gamma ( V \cdot n)$, which results in
\begin{align*}
   \fracp{w'}{n} &= (V \cdot n) \l[\fracp{w_N}{n} - \fracp{^2 w}{n^2} \r]  + \gradient w \cdot \gradient_\Gamma ( V \cdot n) \\
   &= (V \cdot n) \l[\fracp{w_N}{n} - \fracp{^2 w}{n^2} \r]  + \gradient_\Gamma w \cdot \gradient_\Gamma ( V \cdot n),
\end{align*}
where the last transformation directly results from $\gradient w \cdot n = 0$ being inserted into the definition of the tangential gradient. This expression is still not entirely in Hadamard form, but we will later use the tangential Green's formula to conclude this transformation.

\section{Shape derivative of the Navier-Stokes equations in strong and variational form}\label{sec:NSSD}
Before adjoint calculus can be used in Section \ref{section:adjoint_calculus} to finalize the Hadamard form, one first needs to establish the corresponding forward problem. Therefore, we now consider the linearization of the compressible Navier--Stokes equations with respect to a variation of the domain. As mentioned above, we distinguish between the Navier--Stokes equations in pointwise and in variational form. Unsurprisingly, both versions of the forward problem lead to distinct linearizations and this section will be used to discuss the respective differences.

\begin{thm}[Shape derivative of the pointwise Navier-Stokes equations]
\label{thm:sd_pointwise_ns}
The local shape derivative $\u'$ of the Navier--Stokes equations in strong form \eqref{eq:ns_pointwise} is given as the solution of
\begin{equation}\label{eq:NSDSPointwise}
 0 = \gradient \cdot \l( \Fcu\u' - \Fvu\u' - \Fvgradu \gradient\u' \r) \text{ in } \Omega,
\end{equation}
where $\F^c_\u := \fracp{\F^c}{\u}, \F^v_\u := \fracp{\F^v}{\u}$ and $\F^v_{\gradient \u} := \fracp{\F^v}{\gradient \u}$.
\begin{prf}
Applying the local shape derivative to both sides of equation \eqref{eq:ns_pointwise} results in
\begin{align*}
0 &= (\gradient \cdot \l( \Fc - \F^v(\u, \gradient\u) \r))' \\
&= \gradient \cdot \l( \Fcu\u' - \Fvu\u' - \Fvgradu\gradient\u' \r).
\end{align*}
\end{prf}
\end{thm}

\begin{thm}[Shape derivative of the variational Navier-Stokes equations]
\label{thm:sd_variational_ns}
The shape derivative of the variational form of the Navier-Stokes equations \eqref{eq:ns_variational} is given by the problem: Find $\u' \in \mathcal{H}$, such that
\begin{equation}
\begin{aligned}
0=& -\so{\u'}{\l[\Fcu - \Fvu \r]\transpose \gradient \v} - \sgw{\Vn \l[\Fc - \Fv\r] }{\gradient \v}  \\
  &- \so{\u'}{\gradient \cdot \l[\l(\Fvgradu\r)\transpose \gradient \v\r]} + \sg{\u'}{n \cdot \l[\l(\Fvgradu\r)\transpose \gradient \v\r]}\\
  &+ \sgf{\u'}{\l[n \cdot \l(\Fcu - \Fvu\r) \r]\transpose \v} - \sgf{\gradient \u'}{\l(n \cdot \Fvgradu\r) \transpose \v} \\
  &+ \sgw{n \cdot \l(\Fc - \Fv\r)'}{\v}  + \int_{\Gamma_W} \Vn \gradient \cdot {\l(\l[\Fc - \Fv\r] \cdot \v \r)} \dd s 
  \quad\quad \forall \v \in \mathcal{H}
\end{aligned}
\end{equation}
\begin{prf}
Let $\mathcal{H}_t$ be the Hilbert-Space corresponding to $\mathcal{H}$ but defined on the perturbed domain $\Omega_t$ rather than $\Omega$. A shape differentiation of the variational form \eqref{eq:ns_variational} results in
\begin{equation}
\begin{aligned}
0=& \ddt{} -\so{\Fc - \Fv}{\gradient \v} + \ddt{}  \sg{n \cdot \l(\Fc - \Fv\r)}{\v} \\
 =& -\so{\l(\Fc - \Fv\r)'}{\gradient \v} -\so{\Fc - \Fv}{\gradient \v'} -\sg{\Vn \l[\Fc - \Fv\r]}{\gradient \v} \\
  & +\sgf{\l[n \cdot \l(\Fc - \Fv\r)\r]'}{\v} + \sgf{n \cdot \l(\Fc - \Fv\r)}{\v'} \\
  & +\sgw{n \cdot \l(\Fc - \Fv\r)'}{\v} + \sgw{n \cdot \l(\Fc - \Fv\r)}{\v'} \\ 
  & + \int_{\Gamma_W} \Vn \gradient \cdot {\l(\l[\Fc - \Fv\r] \cdot \v \r)} \dd s
\quad \quad \forall \v \in \mathcal{H}, 
\end{aligned}
\end{equation}
where we used equation \eqref{eq:sd_volume} and \eqref{eq:sd_boundary} for the volume and the farfield integrals and Theorem \ref{thm:preliminary_sd} for the wall boundary integral.
The terms containing $\v'$ vanish due to the forward equation \eqref{eq:ns_variational} being satisfied. Using the product rule on $\l[n \cdot \l(\Fc - \Fv\r) \r]'$ and the chain rule on $\l(\F^c(\u)\r)'$ and $\l(\F^v(\u, \gradient \u) \r)'$ at the farfield boundary and in the volume leads to 
\begin{equation*}
\begin{aligned}
0=& -\so{\Fcu \u' - \Fvu \u' - \Fvgradu \gradient \u'}{\gradient \v} -\sg{\Vn \l[\Fc - \Fv\r]}{\gradient \v} \\
  & +\sgf{n \cdot \l[ \Fcu \u' - \Fvu \u' - \Fvgradu \gradient \u' \r]}{\v} + \sgf{n' \cdot \l[ \l(\Fc - \Fv\r)\r]}{\v}  \\
  & +\sgw{n \cdot \l(\Fc - \Fv\r)'}{\v} + \int_{\Gamma_W} \Vn \gradient \cdot {\l(\l[\Fc - \Fv\r] \cdot \v \r)} \dd s
\quad \quad \forall \v \in \mathcal{H}.
\end{aligned}
\end{equation*}
Since $V = 0$ is fulfilled in the neighbourhood of the farfield boundary $\Gamma \setminus \Gamma_W$ the local shape derivative of normal vector $n'$ vanishes and it remains 
\begin{equation*}
\begin{aligned}
0=& -\so{\Fcu \u' - \Fvu \u' - \Fvgradu \gradient \u'}{\gradient \v} -\sgw{\Vn \l[\Fc - \Fv\r]}{\gradient \v} \\
  & +\sgf{n \cdot \l[ \Fcu \u' - \Fvu \u' - \Fvgradu \gradient \u' \r]}{\v} \\
  & +\sgw{n \cdot \l(\Fc - \Fv\r)'}{\v} + \int_{\Gamma_W} \Vn \gradient \cdot {\l(\l[\Fc - \Fv\r] \cdot \v \r)} \dd s
\quad \quad \forall \v \in \mathcal{H}.
\end{aligned}
\end{equation*}
We shift $n$, $\Fcu$, $\Fvu$ and $\Fvgradu$ to the other side of the products and integrate the volume term containing $\gradient \u'$ by parts to obtain the stated expression. 
\end{prf}
\end{thm}

\section{Adjoint calculus}
\label{section:adjoint_calculus}
We now recall adjoint calculus to reformulate a shape optimization problem as documented in \cite{giles1997adjoint}, \cite{jameson1995optimum} or \cite{jameson1998optimum}. In our case the cost function $J$ to be shape optimized is the drag or lift coefficient
\begin{equation*}
J = J(\u, S)
,
\end{equation*}
which depends on a function $S$ describing the shape and the flow solution $\u$ of the governing equation 
\begin{equation}\label{eq:AdjointGeneralStrong}
N(\u, S) = 0
.
\end{equation}
Since $\u$ depends, through the governing equation, on the shape function $S$, a variation of the shape $\delta S$ leads to the following variation of the cost function 
\begin{equation}
\label{eq:pre_deltaJ}
\delta J = \fracp{J}{\u} \delta \u + \fracp{J}{S} \delta S
.
\end{equation}
Therefore, to compute the variation $\delta J$, one needs to know the sensitivity of the flow solution $\delta \u$ for each degree of freedom within the shape deformation. To calculate this variation for each such parameter, a flow solution has to be computed. This prohibitive numerical effort of multiple flow computations can be omitted if it is possible to eliminate the variation $\delta \u$. 
For this purpose, we look at the variation of the governing equation
\begin{equation}
\label{eq:deltaN}
\delta N = \fracp{N}{\u} \delta \u + \fracp{N}{S} \delta S = 0
,
\end{equation}
which provides another equation determining the variation $\delta \u$. Because the variation $\delta N$ equals zero, it can be multiplied by a Lagrange multiplier $\z$ and then be subtracted from the variation of the cost function:
\begin{equation*}
\delta J = \delta J - \z \transpose \delta N
.
\end{equation*}
Inserting equations \eqref{eq:pre_deltaJ} and \eqref{eq:deltaN} for the variations $\delta J$ and $\delta N$ into this equation yields
\begin{equation*}
\begin{aligned}
\delta J &= \fracp{J}{\u} \delta \u + \fracp{J}{S} \delta S - \z \transpose \l(\fracp{N}{\u} \delta \u + \fracp{N}{S} \delta S \r) 
\\
&= \l(\fracp{J}{\u} - \z \transpose \fracp{N}{\u} \r) \delta \u + \l(\fracp{J}{S} - \z \transpose \fracp{N}{S}  \r) \delta S 
.
\end{aligned}
\end{equation*}
The first term, containing $\delta \u$, is then eliminated if $\z$ satisfies the adjoint equation
\begin{equation}
\label{eq:adjoint_equation}
\fracp{J }{\u} - \z \transpose \fracp{N}{\u} = 0
.
\end{equation}
Therewith the variation of the cost function becomes
\begin{equation}
\label{eq:deltaJ}
\delta J = \l(\fracp{J}{S} - \z \transpose \fracp{N}{S}  \r) \delta S,
\end{equation} 
which can be computed without multiple primal flow solutions.

Although the nature of equation~\eqref{eq:AdjointGeneralStrong} makes the above motivation more reminiscent of the strong form approach, it is nevertheless presented here to illustrate the methodology. Furthermore, as stated before, we are also considering the strong form for our shape optimization problem so we have a procedure to gauge against. Thus, we are going to consider the shape derivative of the objective function as given by Theorem~\ref{thm:sd_aerodynamic_targetfunction}, but under the assumption of a state equation in weak form, equation~\eqref{eq:ns_variational}.

\subsection{Variational formulation of the continuous adjoint problem}
We will use this section to derive the variational formulation of the adjoint problem. For more details on the variational approach also see~\cite{KunishVariationalSD} and the respective integral transforms are covered in more depth in~\cite{Har08b}.

\begin{thm}[Variational form of the adjoint Navier--Stokes equations]
\label{thm:variational_adjoint_ns}
The variational formulation of the adjoint Navier--Stokes equations is given by find $\z\in \mathcal{H}$ such that
\begin{multline}
\label{eq:variational_adjoint_ns}
 - \so{\w}{(\F_\u^c - \F_\u^v) \transpose \gradient \z}  
  -\so{\w}{\gradient \cdot \l( \l(\F_{\gradient \u}^v \r) \transpose \gradient \z\r)}
   + \sg{\w}{n \cdot \l( \l(\F_{\gradient \u}^v \r) \transpose \gradient \z \r)}\\
  + \sg{\w}{\l(n \cdot \l(\F_\u^c - \F_\u^v \r) \r)\transpose \z}
   - \sg{\gradient \w}{ \l(n \cdot \F_{\gradient \u}^v\r) \transpose \z}
= J'[\u](\w) \quad \forall \w \in \mathcal{H},
\end{multline}
where the linearization of the cost function in case of drag or lift coefficient is given by
\begin{equation*}
J'[\u](\w) = \sgw{\tfrac{1}{C_\infty}  \l( p_\u n - \tau_\u n \r) \cdot \psi  }{\w } - \sgw{\tfrac{1}{C_\infty}  \l( \tau_{\gradient \u} n \r) \cdot \psi  }{\gradient \w}.
\end{equation*}
\end{thm}
\begin{prf}
Following the outline given by~\cite{KunishVariationalSD}, the linearization in direction $\w \in \mathcal{H}$ of the Navier--Stokes equations in variational form~\eqref{eq:ns_variational} is given by
\begin{align*}
&\langle F'(\u, \Omega) \w, \z\rangle_{\mathcal{H}^* \times \mathcal{H}}\\
 =& -\so{\fracp{\F^c}{\u} \w - \fracp{\F^v}{\u} \w - \fracp{\F^v}{\gradient \u} \gradient \w}{\gradient \z} + \sg{n \cdot \l( \fracp{\F^c}{\u} \w - \fracp{\F^v}{\u} \w - \fracp{\F^v}{\gradient \u} \gradient \w\r) }{\z} \quad \forall \z \in \mathcal{H}
\end{align*}
The adjoint equation in variational form is then given by the problem of finding $\z \in \mathcal{H}$ such that
\begin{equation}
\label{eq:linearized_adjoint}
\langle F'(\u, \Omega) \w, \z\rangle_{\mathcal{H}^* \times \mathcal{H}} = J'[\u](\w) \quad \forall \w \in \mathcal{H}.
\end{equation}
Thus, we have
\begin{align*}
 &- \so{(\F_\u^c - \F_\u^v) \w - \F_{\gradient \u}^v \gradient \w}{\gradient \z} 
 + \sg{n \cdot ((\F_\u^c - \F_\u^v) \w  - \F_{\gradient \u}^v \gradient \w )}{\z}\\
=& - \so{\w}{(\F_\u^c - \F_\u^v) \transpose \gradient \z} 
   + \so{\gradient \w}{ \l(\F_{\gradient \u}^v \r) \transpose \gradient \z} + \sg{\w}{\l(n \cdot \l(\F_\u^c - \F_\u^v \r) \r)\transpose \z}
   - \sg{\gradient \w}{ \l(n \cdot \F_{\gradient \u}^v\r) \transpose \z} \\
=& J'[\u](\w) \quad \forall \w \in \mathcal{H}.
\end{align*}
Integration by parts resolves the remaining gradient $\gradient \w$ in the second volume integral and the variational formulation of the adjoint Navier--Stokes equations becomes the desired expression.
\end{prf}
\begin{kor}
Choosing $\w$ in equation \eqref{eq:variational_adjoint_ns} with appropriate compact support in either $\Omega$ or on $\Gamma_W$ and $\Gamma\setminus\Gamma_W$, one can see that
\begin{align} 
\label{eq:volume_condition_z}
 - \so{\w}{(\F_\u^c - \F_\u^v) \transpose \gradient \z} - \so{\w}{\gradient \cdot \l( \l(\F_{\gradient \u}^v \r) \transpose \gradient \z\r)} =&  0 \quad \forall \w \in \mathcal{H}_0(\Omega)
\end{align}
for the volume. For a test function with compact support on $\Gamma_W$ we see that
\begin{align}\begin{aligned}
 \label{eq:AdjWingBoundary}
   &\sgw{\w}{n \cdot \l( \l(\F_{\gradient \u}^v \r) \transpose \gradient \z \r)} + \l(\w, \l(n \cdot \l(\F_\u^c - \F_\u^v \r) \r)\transpose \z\r)_{\Gamma_W} - \l(\gradient \w, \l(n \cdot \F_{\gradient \u}^v\r) \transpose \z\r)_{\Gamma_W} \\
  =& \sgw{\tfrac{1}{C_\infty}  \l( p_\u n - \tau_\u n \r) \cdot \psi  }{\w } - \sgw{\tfrac{1}{C_\infty}  \l( \tau_{\gradient \u} n \r) \cdot \psi  }{\gradient \w} \quad \forall \w \in \mathcal{H} \cap \mathcal{H}_0(\Gamma_W)
\end{aligned}\end{align}
and finally using the same argumentation on all remaining boundaries $\Gamma\setminus\Gamma_W$
\begin{multline}
\label{eq:boundary_condition_z}
   \sgf{\w}{n \cdot \l( \l(\F_{\gradient \u}^v \r) \transpose \gradient \z \r)} + \sgf{\w}{\l(n \cdot \l(\F_\u^c - \F_\u^v \r) \r)\transpose \z} - \sgf{\gradient \w}{\l(n \cdot \F_{\gradient \u}^v\r) \transpose \z} \\
   = 0 \quad \forall \w \in \mathcal{H} \cap \mathcal{H}_0(\Gamma\setminus\Gamma_W). 
\end{multline}
\end{kor}

\subsection{Application of adjoint equation to the shape derivative of the Navier--Stokes equations}
Recalling our goal of eliminating the local shape derivatives $p'$ and $\tau'$ in equation~\eqref{eq:PrelimGradient}, we will now derive two intermediate relationships between the adjoint equation and the respective linearizations of the Navier--Stokes equations. One stems from a consideration of a pointwise linearization while the other stems from the same procedure applied to the linearization of the weak form. The resulting two different intermediate expressions~\eqref{eq:sd_ns_pointwise} and~\eqref{eq:sd_ns_variational} will then be transformed further in Section~\ref{subsec:BoundaryDerivation}. There, the corresponding variations in terms of the pressure-based variables will be made explicit. The variation of the boundary conditions will be taken into consideration in Section~\ref{sec:BoundarySubstraction}.

We begin by considering the pointwise problem. Multiplying the shape derivative of pointwise Navier-Stokes equations from Theorem \ref{thm:sd_pointwise_ns} by a testfunction $\v$ and integrating over the domain $\Omega$ gives us 
\begin{equation*}
 0 = \so{\gradient \cdot \l( \Fcu\u' - \Fvu\u' - \Fvgradu \gradient\u' \r)}{\v} 
 \quad \forall \v \in \mathcal{H}.
\end{equation*}
Integration by parts yields
\begin{equation*}
\begin{aligned}
 0 =& -\so{\l( \Fcu\u' - \Fvu\u' - \Fvgradu \gradient\u' \r)}{\gradient \v} \\
    & +\sgf{n \cdot \l( \Fcu\u' - \Fvu\u' - \Fvgradu \gradient\u' \r)}{\v} \\
    & +\sgw{n \cdot \l( \Fcu\u' - \Fvu\u' - \Fvgradu \gradient\u' \r)}{\v}
 \quad \forall \v \in \mathcal{H}.
\end{aligned}
\end{equation*}
Shifting $n$, $\Fcu$, $\Fvu$ and $\Fvgradu$ to the other side of the products results in
\begin{equation*}
\begin{aligned}
 0 =& -\so{\u'}{ \l[ \Fcu - \Fvu \r]\transpose \gradient \v} + \so{\gradient \u'}{ \l( \Fvgradu \r) \transpose \gradient \v} \\
    & +\sgf{\u'}{ \l[ n \cdot \l( \Fcu - \Fvu \r) \r] \transpose \v} -\sgf{\gradient \u'}{ \l[n \cdot \l(\Fvgradu \r) \r]\transpose \v} \\
    & +\sgw{n \cdot \l( \Fcu\u' - \Fvu\u' - \Fvgradu \gradient\u' \r)}{\v}
 \quad \forall \v \in \mathcal{H}.
\end{aligned}
\end{equation*}
Integration by parts in the 2nd volume integral and applying the chain rule backwards to the wall integral leads to
\begin{equation*}
\begin{aligned}
 0 =& -\so{\u'}{ \l[ \Fcu - \Fvu \r]\transpose \gradient \v} \\
    & -\so{\u'}{\gradient \cdot \l[\l( \Fvgradu \r) \transpose \gradient \v\r]} + \sg{\u'}{n \cdot \l[ \l( \Fvgradu \r) \transpose \gradient \v \r]} \\
    & +\sgf{\u'}{ \l[ n \cdot \l( \Fcu - \Fvu \r) \r] \transpose \v} -\sgf{\gradient \u'}{ \l[n \cdot \l(\Fvgradu \r) \r]\transpose \v} \\
    & +\sgw{n \cdot \l(\Fc - \Fv\r)'}{\v}
 \quad \forall \v \in \mathcal{H}.
\end{aligned}
\end{equation*}
Using adjoint conditions \eqref{eq:volume_condition_z} and \eqref{eq:boundary_condition_z} and changing the name of the dependent variable from $\v$ to $\z$, one obtains
\begin{equation}
\label{eq:sd_ns_pointwise}
 0 = \sgw{\u'}{n \cdot \l[ \l( \Fvgradu \r) \transpose \gradient \z \r]} + \sgw{n \cdot \l(\Fc - \Fv\r)'}{\z}
 \quad \forall \z \in \mathcal{H}.
\end{equation}
Using equations \eqref{eq:volume_condition_z} and \eqref{eq:boundary_condition_z}, which stem from a weak form adjoint equation, within the strong form linearization here might appear counter intuitive at first glance. However, it should be noted that a pointwise interpretation of those does not effect the above equation.

Contrary to the above preliminary result stemming from the strong form of the Navier--Stokes equations, we next follow the same process, now considering the interaction with the adjoint equation and the variational Navier--Stokes equations from Theorem \ref{thm:sd_variational_ns}, which results in
\begin{equation*}
\begin{aligned}
0=&\sgw{\u'}{n \cdot \l[\l(\Fvgradu\r)\transpose \gradient \z\r]} + \sgw{n \cdot \l(\Fc - \Fv\r)'}{\z} \\
  &+ \int_{\Gamma_W} \Vn \gradient \cdot {\l(\l[\Fc - \Fv\r] \cdot \z \r)} \dd s - \sgw{\Vn \l[\Fc - \Fv\r] }{\gradient \z} 
  \quad\quad \forall \z \in \mathcal{H}.
\end{aligned}
\end{equation*}
We apply the product rule to the divergence and get
\begin{equation}
\begin{aligned}
\label{eq:sd_ns_variational}
0=&\sgw{\u'}{n \cdot \l[\l(\Fvgradu\r)\transpose \gradient \z\r]} + \sgw{n \cdot \l(\Fc - \Fv\r)'}{\z} \\
  &+ \extrabox{\sgw{\Vn \gradient \cdot {\l[\Fc - \Fv\r]}}{\z} }
  \quad\quad \forall \z \in \mathcal{H}.
\end{aligned}
\end{equation}
Comparing equation \eqref{eq:sd_ns_pointwise} and \eqref{eq:sd_ns_variational} we see that for the variational form of the Navier-Stokes equations, there is one extra term (framed), which vanishes if the Navier--Stokes equations are fulfilled pointwise. In the following, we can therefore avoid a fork and always use equation \eqref{eq:sd_ns_variational}, keeping in mind that all framed terms only occur in the variational approach. 
Also, expressing the variation $\u'$ in terms of the non-conserved variables works likewise, irrespective of the underlying form of the state equation.

\subsection{Transformation to non-conserved variables}\label{subsec:BoundaryDerivation}

In this subsection we insert in particular the no-slip condition into the shape derivative of the Navier-Stokes equations \eqref{eq:sd_ns_variational}. Therefore we first state the so called homogeneity tensor $G = \F^v_{\gradient \u}$ at the no-slip wall. We also use $v=0$ to get the representations of $\u', (\F^c)'$ and $(\F^v)'$ at the wall boundary. With these terms the local shape derivative $\u'$ is then reformulated, such that we form the respective local shape derivative of the pressure $p'$ and the viscous stress tensor $\tau'$, which will later eliminate their counterparts in the preliminary shape derivative of the cost function, given by equation~\eqref{eq:PrelimGradient}.

As in~\cite{Har08b}, the homogeneity tensor is given by 
\begin{equation*}
G = \l[G^{ij}_{kl}\r]^{ij}_{kl} = \fracp{\l(f^v_k\r)_i}{\fracp{\u_j}{x_l}}
,
\end{equation*}
where $(f^v_k)_i$ denotes the $i$-th component of the $k$-th viscous flux vector and $\fracp{\u_j}{x_l}$ denotes the derivative of the $j$-th component of the vector of conserved variables with respect to $x_l$. 
Therefore in two dimensions the indices fulfill $i,j \in \{1, \ldots, 4\}$ and $k, l \in \{1, 2\}$. 
The homogeneity tensor $G$ at the no-slip wall is given by
\begin{equation*}
\begin{aligned}
G_{11} &= \frac{\mu}{\rho} 
   \begin{pmatrix}
      0 &  0  &  0  &  0\\ 
      0 & \tfrac{4}{3} & 0 & 0\\
      0 & 0 & 1 & 0\\
      -\tfrac{\gamma}{Pr}E & 0 & 0 & \tfrac{\gamma}{Pr}\\
   \end{pmatrix}, 
G_{12} = \frac{\mu}{\rho} 
   \begin{pmatrix}
      0  &  0  &  0  &  0\\
      0  &  0  & -\tfrac{2}{3} & 0\\
      0 & 1 & 0 & 0\\
      0  &  0  &  0  &  0
   \end{pmatrix}, \\
G_{21} &= \frac{\mu}{\rho} 
   \begin{pmatrix}
      0  &  0  &  0 &  0\\
      0  &  0  &  1 &  0\\
      0  &  -\tfrac{2}{3} & 0 & 0\\
      0  &  0  &  0  &  0
   \end{pmatrix}, 
G_{22} = \frac{\mu}{\rho} 
   \begin{pmatrix}
      0 &  0  &  0  &  0\\ 
      0 &  1  &  0  & 0\\
      0 &  0  & \tfrac{4}{3} & 0\\
      -\tfrac{\gamma}{Pr}E & 0 & 0 & \tfrac{\gamma}{Pr}\\
   \end{pmatrix}.
\end{aligned}
\end{equation*}
See~\cite{Har08b} for a more detailed discourse on the general homogeneity tensor. 
%
The shape derivative of the vector of conserved variables $\u'$ reduces with $(\rho v_i)' = \rho v_i' + \rho' v_i = \rho v_i'$ at the no-slip wall to
\begin{equation*}
\u' = \begin{pmatrix}
\rho \\ \rho v_1 \\ \rho v_2 \\ \rho E
\end{pmatrix} ' = 
\begin{pmatrix}
\rho' \\ \rho v'_1 \\ \rho v'_2 \\ (\rho E)' 
\end{pmatrix}
 \quad \on \Gamma_W
.
\end{equation*}
%
Because $(\rho v_i v_j)' = (\rho v_i)' v_j + \rho v_i v'_j = 0$ holds due the no-slip condition, the shape derivative of the convective flux $(\F^c)'$ reduces to 
\begin{equation}
\label{eq:(F^c)'}
(\F^c)' = \begin{pmatrix}
\rho v_1 & \rho v_2 \\ 
\rho v_1^2 + p & \rho v_1 v_2 \\ 
\rho v_1 v_2 & \rho v_2^2 + p \\
\rho H v_1 & \rho H v_2
\end{pmatrix} '  = \begin{pmatrix}
\rho v'_1 & \rho v'_2 \\ 
p' & 0 \\ 
0 & p' \\
\rho H v'_1 & \rho H v'_2
\end{pmatrix}
\quad \on \Gamma_W
.
\end{equation}
%
For the shape derivative of the viscous flux $(\F^v)'$ we can use $(\tau_{ij} v_j)' = \tau_{ij} v'_j$ at the no-slip boundary to obtain
\begin{equation}
\label{eq:(F^v)'}
(\F^v)' = \begin{pmatrix}
0 & 0 \\
\tau_{11}' & \tau_{12}' \\ 
\tau_{21}' & \tau_{22}' \\ 
\sum_j \tau_{1j} v'_j + \kappa \fracp{T'}{x_1} & \sum_j \tau_{2j} v'_j + \kappa \fracp{T'}{x_2}
\end{pmatrix}
\quad \on \Gamma_W
.
\end{equation}

These terms are now inserted into the shape derivative of the Navier-Stokes equations for each scalar product of equation \eqref{eq:sd_ns_variational} seperately. 
First we consider the second scalar product of equation \eqref{eq:sd_ns_variational} and insert the representations \eqref{eq:(F^c)'} and \eqref{eq:(F^v)'} of $(\F^c)'$ and $(\F^v)'$ into this integral:
\begin{align*}
&\sgw{ n \cdot \l( \l( \F^c \r)' - \l( \F^v \r)' \r)}{\z}
=&\s{\l[
\begin{pmatrix}
\rho v'_1 & \rho v'_2 \\ 
p' & 0 \\ 
0 & p' \\
\rho H v'_1 & \rho H v'_2
\end{pmatrix}
- \begin{pmatrix}
0 & 0 \\
\tau_{11}' & \tau_{12}' \\ 
\tau_{21}' & \tau_{22}' \\ 
\sum_j\tau_{1j} v'_j + \kappa \fracp{T'}{x_1} & \sum_j\tau_{2j} v'_j + \kappa \fracp{T'}{x_2}
\end{pmatrix}
\r] n}{\z}{\Gamma_W}
.
\end{align*}
From the first line we read $\rho (v' \cdot n ) \z_1$, from the second and third line $(p'n - \tau'n) \cdot \z_{2,3}$ and the last line gives $\l(\rho H v' - \l( \tau v' + \kappa \gradient T' \r)\r) \cdot n \z_4$. Altogether this yields 
\begin{align}
\begin{aligned}
\label{item:2.scalarproduct}
&\sgw{ n \cdot \l( \l( \F^c \r)' - \l( \F^v \r)' \r)}{\z}
=&\int_{\Gamma_W} (p'n - \tau' n) \cdot \z_{2,3} + v' \cdot \l( \rho n \z_1 + ( \rho H n - \tau n) \z_4 \r) - \gradient T' \cdot n \kappa \z_4 \dd s.
\end{aligned}
\end{align}
%
Next we use the homogeneity tensor to rewrite the first scalar product of equation \eqref{eq:sd_ns_variational}:
\begin{equation*}
\sgw{\u'}{n \cdot \l(\l(\F_{\gradient \u}^v \r)\transpose \gradient \z \r)} 
= \sgw{\u'}{n \cdot \l(G\transpose \gradient \z \r)} = \int_{\Gamma_W} \sum_{k,l,i,j} \u'_j n_l G^{ij}_{kl} \fracp{\z_i}{x_k} \dd s
.
\end{equation*}
As we can see, the first component of $\z$ is multiplied with the first lines of the matrices $G_{kl}$, which are all zero. Therefore there is no contribution for $\z_1$. The non-zero entries in the second and third line of $G$ at the no-slip wall look familiar to the coefficients of $\fracp{v_i}{x_j}$ in the viscous stress tensor. If we define the \definition{adjoint stress tensor} as 
\begin{equation*}
\Sigma := \mu \l(\gradient \z_{2,3} + (\gradient \z_{2,3})\transpose - \frac{2}{3} \l(\gradient \cdot \z_{2,3} \r) I \r)
\end{equation*}
one can actually see that the second and third line yield
\begin{equation*}
\int_{\Gamma_W} \sum_{\substack{i,j=2,3 \\ k,l}} \u'_j n_l G^{ij}_{kl} \fracp{\z_i}{x_k} \dd s
= \int_{\Gamma_W} v' \cdot (n \cdot \Sigma) \dd s
.
\end{equation*}
For the fourth component of $\z$, meaning the fourth lines of $G$, we only get contributions for $k=l$ and $j =1$ or $4$, meaning the first and last entry in the fourth line of the matrices $G_{11}$ and $G_{22}$. Since the $j$-th columns of the matrices $G_{kl}$ are multiplied by the $j$-th component of $\u$, the remaining expression is
\begin{equation*}
\int_{\Gamma_W} \sum_{\substack{i=4 \\ j=1,4 \\ k=l}} \u'_j n_l G^{ij}_{kl} \fracp{\z_i}{x_k} \dd s
= \int_{\Gamma_W} \frac{\mu}{\rho} \frac{\gamma}{Pr} \l(-E\u'_1 + \u'_4 \r) n \cdot \gradient \z_4 \dd s
.
\end{equation*}
With $\u'_1 = \rho'$ and $\u'_4 = (\rho E)' = \rho' E + \rho E'$ this integral is simplified to 
\begin{equation*}
\int_{\Gamma_W} \frac{\mu}{\rho} \frac{\gamma}{Pr} \rho E' n \cdot \gradient \z_4 \dd s = \int_{\Gamma_W} T' \kappa n \cdot \gradient \z_4 \dd s
,
\end{equation*}
where we use $T \kappa = \frac{\mu \gamma}{Pr} \l(E - \frac{1}{2} v^2 \r)$ and the no-slip condition to substitute
\begin{equation*}
T' \kappa = \frac{\mu \gamma}{Pr} \l(E' - \l( \frac{1}{2} v^2 \r)' \r) = \frac{\mu \gamma}{Pr} E'
.
\end{equation*}
In total we obtain for the first scalar product of equation \eqref{eq:sd_ns_variational} 
\begin{equation} \label{item:1.scalarproduct}
\sgw{\u'}{n \cdot \l(\l(\F_{\gradient \u}^v \r)\transpose \gradient \z \r)} 
= \int_{\Gamma_W} v' \cdot (n \cdot \Sigma) \dd s + \int_{\Gamma_W} T' \kappa n \cdot \gradient \z_4 \dd s
.
\end{equation}
Combining the results of \eqref{item:2.scalarproduct} and \eqref{item:1.scalarproduct} the preliminary shape derivative equation~\eqref{eq:sd_ns_variational}, now containing the variation of the boundary condition and the explicit variations of the primal variables, becomes
\begin{multline*}
0=\int_{\Gamma_W} v' \cdot (n \cdot \Sigma) \dd s + \int_{\Gamma_W} T' \kappa n \cdot \gradient \z_4 \dd s \\
+ \int_{\Gamma_W} (p'n - \tau' n) \cdot \z_{2,3} + v' \cdot \l( \rho n \z_1 + ( \rho H n - \tau n) \z_4 \r) - \gradient T' \cdot n \kappa \z_4 \dd s \\
+ \extrabox{\sgw{(V \cdot n) (\gradient \cdot (\F^c - \F^v))}{\z}}.
\end{multline*}

\subsection{Subtraction of the shape derivative of the Navier-Stokes equations from the preliminary shape derivative of the cost function}\label{sec:BoundarySubstraction}
The above equality is now, in accordance with the adjoint approach, subtracted from the preliminary shape derivative of the cost function~\eqref{eq:PrelimGradient}, thereby obtaining a representation that does not contain any local shape derivatives.
\begin{equation}
\label{eq:proof_hadamard_costfunction1}
\begin{aligned}
\dJ =& \frac{1}{C_\infty} \int_{\Gamma_W} {(p' n - \tau' n ) \cdot \psi}^\num{1}  + {(V \cdot n) \div(p\psi - \tau \psi)}^\num{2} \dd s \\
&- \int_{\Gamma_W} v' \cdot (n \cdot \Sigma) \dd s - \int_{\Gamma_W} T' \kappa n \cdot \gradient \z_4 \dd s \\
&- \int_{\Gamma_W} {(p'n - \tau' n) \cdot \z_{2,3}}^\num{3} + v' \cdot \l( \rho n \z_1 + ( \rho H n - \tau n) \z_4 \r) - \gradient T' \cdot n \kappa \z_4 \dd s \\
&- {\extrabox{\sgw{(V \cdot n) (\gradient \cdot (\F^c - \F^v))}{\z}}}^\num{4} 
\end{aligned}
\end{equation}
To ease the following discussions we use framed numbers $\num{1}, \num{2}, \num{3}, \ldots$ to refer to certain terms. Term $\num{3}$ and term $\num{1}$ cancel each other if the adjoint boundary condition $\z_{2,3} = \frac{1}{C_\infty} \psi \on \Gamma_W$ is fulfilled. Furthermore, this boundary condition also implies that the $\z_{2,3}$-component of expression \num{4}   cancels out term \num{2}, because due to $\fracp{\rho v_i v_j}{x_k} = \fracp{\rho v_i}{x_k} v_j + \rho v_i \fracp{v_j}{x_k} = 0$ at the no-slip boundary, the equation 
\begin{equation}
\label{eq:NS_23}
\l(\gradient \cdot (\F^c - \F^v)\r)_{2,3} = \gradient \cdot \l[(\rho v_i v_j)_{ij} + pI - \tau \r] = \div(pI - \tau)
\end{equation}
and therewith $\l(\gradient \cdot (\F^c - \F^v)\r)_{2,3} \cdot \z_{2,3} = \frac{1}{C_\infty} \div( p \psi - \tau \psi)$ holds.

When considering the pointwise equations the scalar product $\num{4}$ is not present and one would expect $\num{2}$ to remain untouched. In this situation, however, equation \eqref{eq:NS_23}, being the pointwise conservation of momentum, equals zero and therefore term $\num{2}$ vanishes anyway. Summarizing the above, the first and last component of $\num{4}$ are the only differences between the shape derivatives of the two approaches. We mark them by $\num{4_{1,4}}$.
\begin{equation}
\label{eq:proof_hadamard_costfunction2}
\begin{aligned}
\dJ =&- \int_{\Gamma_W} v' \cdot (n \cdot \Sigma) \dd s - \int_{\Gamma_W} T' \kappa n \cdot \gradient \z_4 \dd s \\
&- \int_{\Gamma_W} v' \cdot \l( \rho n \z_1 + ( \rho H n - \tau n) \z_4 \r) - \gradient T' \cdot n \kappa \z_4 \dd s \\
&- {\extrabox{\sgw{(V \cdot n) (\gradient \cdot (\F^c - \F^v))_{1,4}}{\z_{1,4}}}}^\num{4_{1,4}} 
\end{aligned}
\end{equation}
To achieve independence of the remaining local shape derivatives $v'$ and $T'$, we use the variation of the boundary conditions, as described in Section \ref{section:sd_of_boundary_conditions}, for the no-slip, adiabatic and isothermal boundary condition:
\begin{equation*}
\begin{aligned}
v' &= - (V \cdot n) \fracp{v}{n} &&\quad \on \Gamma_W, \\
T' &= (V \cdot n) \fracp{T_W - T}{n} &&\quad \on \Gamma_{iso}, \\
\gradient T' \cdot n &= - (V \cdot n) \fracp{^2 T}{n^2} + \gradient T \cdot \gradient_\Gamma( V \cdot n)  &&\quad \on \Gamma_{adia}
.
\end{aligned}
\end{equation*}
Inserting these conditions on the respective parts of the wall, equation \eqref{eq:proof_hadamard_costfunction2} becomes 
\begin{equation}
\label{eq:proof_hadamard_costfunction3}
\begin{aligned}
\dJ =& \int_{\Gamma_W} (V \cdot n) \fracp{v}{n} \cdot (n \cdot \Sigma) \dd s \\
&- \int_{\Gamma_{iso}} (V \cdot n) \fracp{T_W - T}{n} \kappa n \cdot \gradient \z_4 \dd s - \int_{\Gamma_{adia}} T' \kappa n \cdot \gradient \z_4 \dd s^\num{5}  \\
&+ \int_{\Gamma_W} (V \cdot n) \fracp{v}{n} \cdot \l( \rho n \z_1 + ( \rho H n - \tau n) \z_4 \r) \dd s \\
&+ \int_{\Gamma_{iso}} \gradient T' \cdot n \kappa \z_4 \dd s^\num{6} - \int_{\Gamma_{adia}} (V \cdot n) \fracp{^2 T}{n^2} \kappa \z_4 \dd s + \int_{\Gamma_{adia}} \gradient T \cdot \gradient_\Gamma ( V \cdot n) \kappa \z_4\dd s^\num{7} \\
&- {\extrabox{\sgw{(V \cdot n) (\gradient \cdot (\F^c - \F^v)_{1,4})}{\z_{1,4}}}}^\num{4_{1,4}} 
.
\end{aligned}
\end{equation}
Application of the tangential Green's formula \eqref{eq:greens_formula} to integral $\num{7}$ leads to
\begin{equation*}
\begin{aligned}
\int_{\Gamma_{adia}} \gradient T \cdot \gradient_\Gamma (V \cdot n) \kappa \z_4 \dd s  
&=\int_{\Gamma_{adia}} (V \cdot n) K (\gradient T \cdot n) \kappa \z_4 - (V \cdot n) \div_\Gamma(\gradient T \kappa \z_4) \dd s \\
&= - \int_{\Gamma_{adia}} (V \cdot n) \div_\Gamma(\gradient T \kappa \z_4) \dd s
,
\end{aligned}
\end{equation*}
where the adiabatic wall condition $\gradient T \cdot n = 0$ was used in the second line.

Integrals of equation \eqref{eq:proof_hadamard_costfunction2} containing the remaining local shape derivatives $T'$ at the adiabatic wall $\num{5}$ and $\gradient T'$ at the isothermal wall $\num{6}$ vanish if the following adjoint boundary conditions are fulfilled
\begin{equation*}
\z_4 = 0, ~ \on \Gamma_{iso} \quad\quad \text{and} \quad\quad \gradient \z_4 \cdot n = 0, ~ \on \Gamma_{adia}
.
\end{equation*}

The last step is to give an expression for the term $\num{4_{1,4}}$. We therefore consider this expression separately for the first and forth component, denoted by $\gradient \cdot (\F^c - \F^v)_{1,4}$. Because the first line of the viscous flux $\F^v$ equals zero anyway and furthermore the no-slip condition holds, the first component simplifies to
\begin{equation*}
\gradient \cdot (\F^c - \F^v)_1 = \gradient \cdot (\rho v_1, \rho v_2) = \rho (\gradient \cdot v)
.
\end{equation*}
For the fourth component we get, again with the no-slip condition, 
\begin{equation*}
\begin{aligned}
\gradient \cdot (\F^c - \F^v)_4 &= \gradient \cdot (\rho H v_1, \rho H v_2) - \gradient \cdot \l(\sum_j \tau_{1j} v_j + \kappa \fracp{T}{x_1}, \sum_j \tau_{2j} v_j + \kappa \fracp{T}{x_2} \r) \\
&= \rho H (\gradient \cdot v) - \sum_{i,j} \tau_{ij} \fracp{v_j}{x_i} - \kappa \Delta T
.
\end{aligned}
\end{equation*}
Altogether the shape derivative of the drag or lift coefficient in Hadamard form becomes
\begin{equation}
\label{eq:hf_targetfunction}
\begin{aligned}
\dJ =& \int_{\Gamma_W} (V \cdot n) \fracp{v}{n} \cdot (n \cdot \Sigma) \dd s \\
&- \int_{\Gamma_{iso}} (V \cdot n) \fracp{T_W - T}{n} \kappa n \cdot \gradient \z_4 \dd s \\
&+ \int_{\Gamma_W} (V \cdot n) \fracp{v}{n} \cdot \l( \rho n \z_1 + ( \rho H n - \tau n) \z_4 \r) \dd s \\
&- \int_{\Gamma_{adia}} (V \cdot n) \l(\fracp{^2 T}{n^2} \kappa \z_4 + \div_\Gamma(\gradient T \kappa \z_4) \r) \dd s \\
&- \extrabox{\int_{\Gamma_W} (V \cdot n) \rho (\gradient \cdot v) \z_1 \dd s}
- \extrabox{\int_{\Gamma_{adia}} (V \cdot n) \l(\rho H (\gradient \cdot v) - \sum_{i,j} \tau_{ij} \fracp{v_j}{x_i} - \kappa \Delta T \r) \z_4 \dd s }
,
\end{aligned}
\end{equation}
where the difference between the variational and the pointwise approach is based on the two framed integrals in the last line, which only appear upon consideration of the variational form of the Navier-Stokes equations.

\section{Numerical results}
We now present numerical results to show the difference between the pointwise and the variational approach in application. The shape derivative of both drag and lift coefficient is implemented in the discontinuous Galerkin solver PADGE developed primarily at the German Aerospace Center (DLR), Braunschweig. This flow solver is already capable to also compute the necessary adjoint solutions, originally implemented for error estimation and $p$-/$h$-refinement. This adjoint solution, together with the primal solution, is directly used to calculate the Hadamard form \eqref{eq:hf_targetfunction} of the drag and lift coefficient. 
Part of the PADGE environment is the ADIGMA MTC3 test case, which is defined as the calculation of flow around a NACA0012 airfoil at Mach $M=0.5$, angle of attack $\alpha = 2.0^\circ$ and a Reynolds number of $Re = 5,000$. 
Because this test case is thoroughly verified, it is a good basis gauge our shape derivatives.  Our grid consists of $1,640$ cells and the profile as the boundary of the mesh is represented by $40$ curved edges of polynomial order four. To obtain a very accurate solution and study the error and behavior for varying polynomial degree $p$ of the Galerkin ansatz, we use a $p$-refinement from degree three to five. The polynomial degree of the adjoint solution was always chosen to be one degree higher, i.e.\ four to six. We verify our two respective shape derivative implementations against finite differences, considering lift and drag separately. The finite difference reference solution is created as follows. Each edge is disturbed by a fourth order polynomial individually, such that the profile stays smooth. A flow solution is then computed for each such perturbation. Because the shape gradient can be evaluated in every support point of these degree-4-polynomials representing the edges of the profile, we calculate four times as many gradient components as edges. 

\begin{figure}[!htb]
\centering
\includegraphics[width=1.0\textwidth]{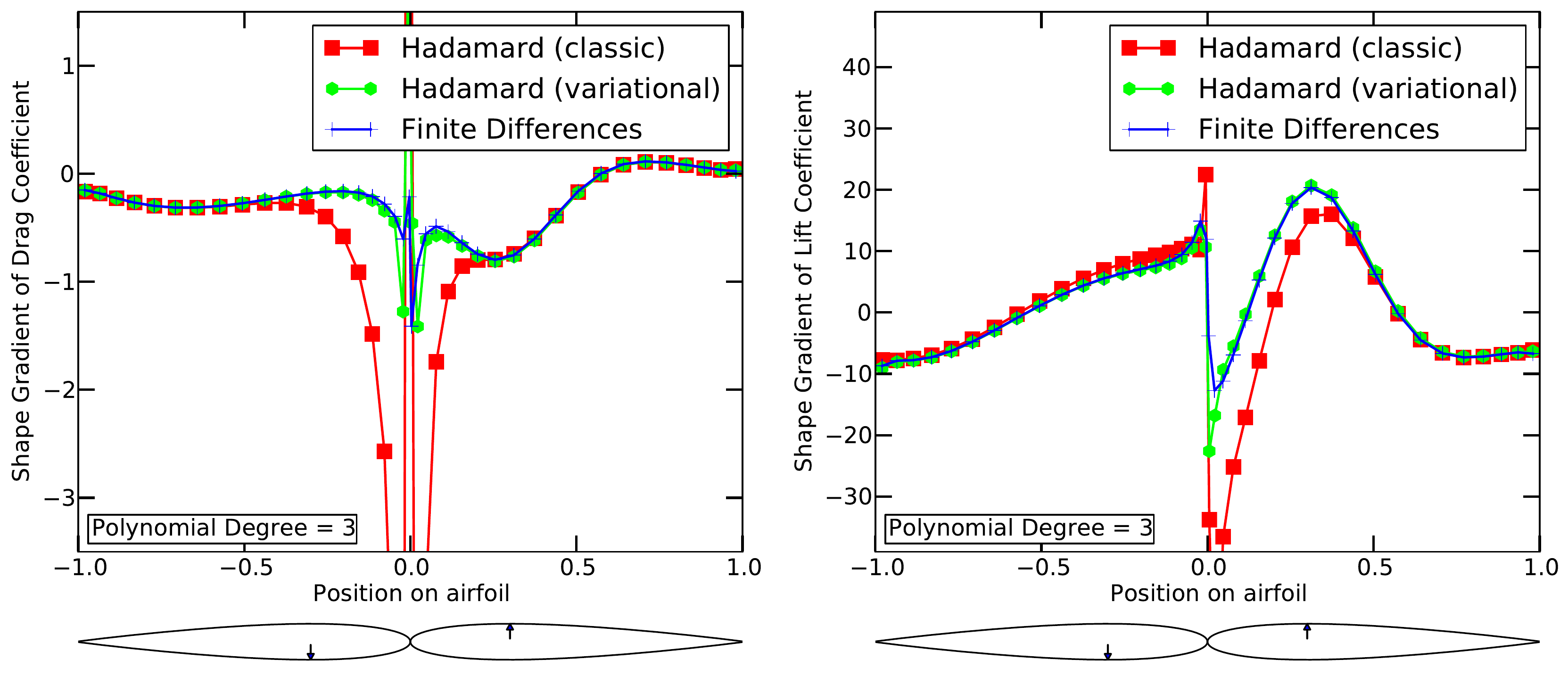}
\caption{Shape Gradient of Drag/Lift Coefficient to polynomial degree $p=3$.}
\label{fi:draglift1}
\end{figure}

With respect to visualizing the results, we always plot the respective data point of the center of each edge. In Figure \ref{fi:draglift1} the finite differences of the drag and lift coefficient are plotted against the shape derivatives in Hadamard form for the variational and the pointwise approach. Both the finite differences and the two Hadamard forms were calculated for a solution made of third order polynomials. One can observe that the shape derivative of the variational approach matches very nicely with the finite differences, whereas the shape derivative of the pointwise approach deviates noticeably.

\begin{figure}[!htb]
\centering
\includegraphics[width=1.0\textwidth]{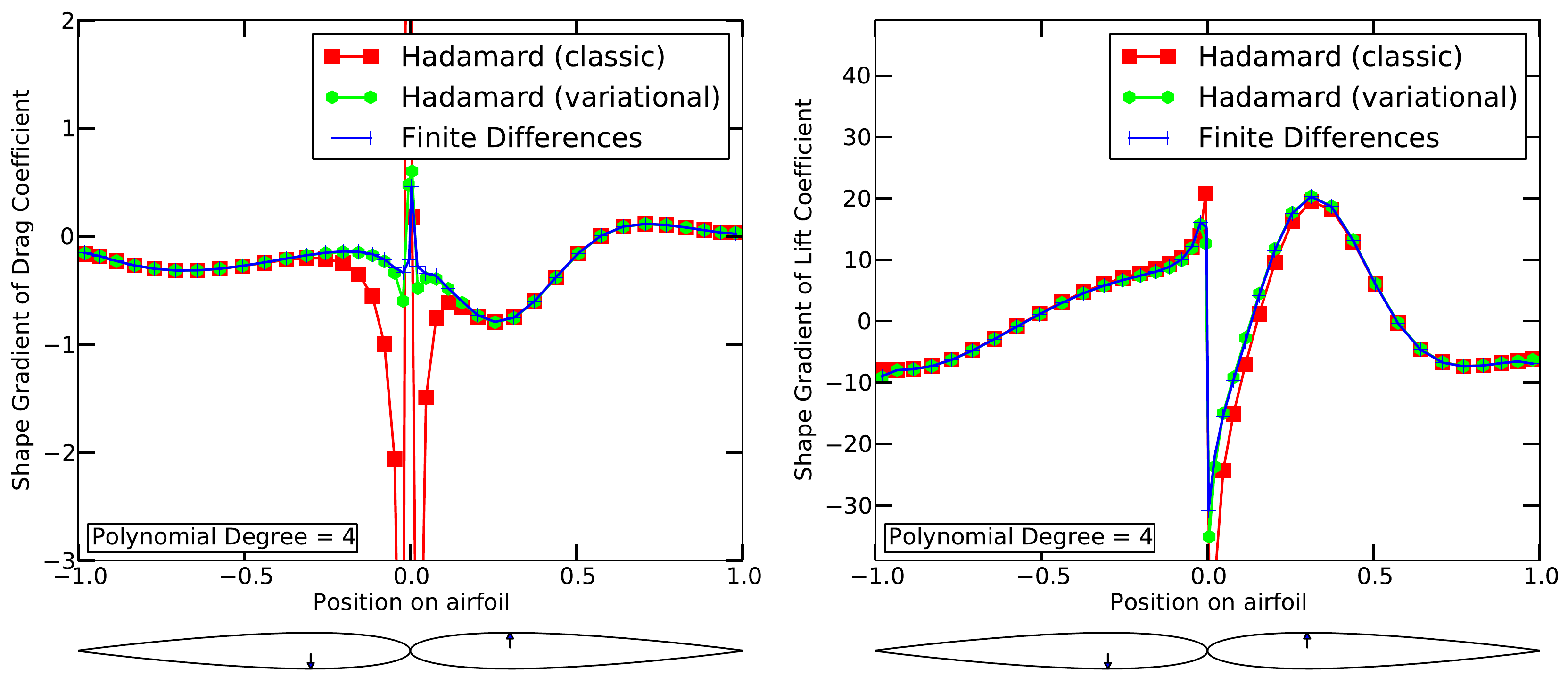}
\caption{Shape Gradient of Drag/Lift Coefficient to polynomial degree $p=4$.}
\label{fi:draglift2}
\end{figure}

As we increase the polynomial degree of the flow solution, see Figure \ref{fi:draglift2}, both Hadamard forms fit the finite differences better overall. However around the nose of the profile and the forward pressure stagnation point, the offset of the Hadamard form of the pointwise approach is still unmissable, most likely due to the magnitude of the gradients of the flow solution. Remarkably, the shape derivative stemming from the variational approach aligns nearly perfectly with the finite difference reference. 

\begin{figure}[!htb]
\centering
\includegraphics[width=1.0\textwidth]{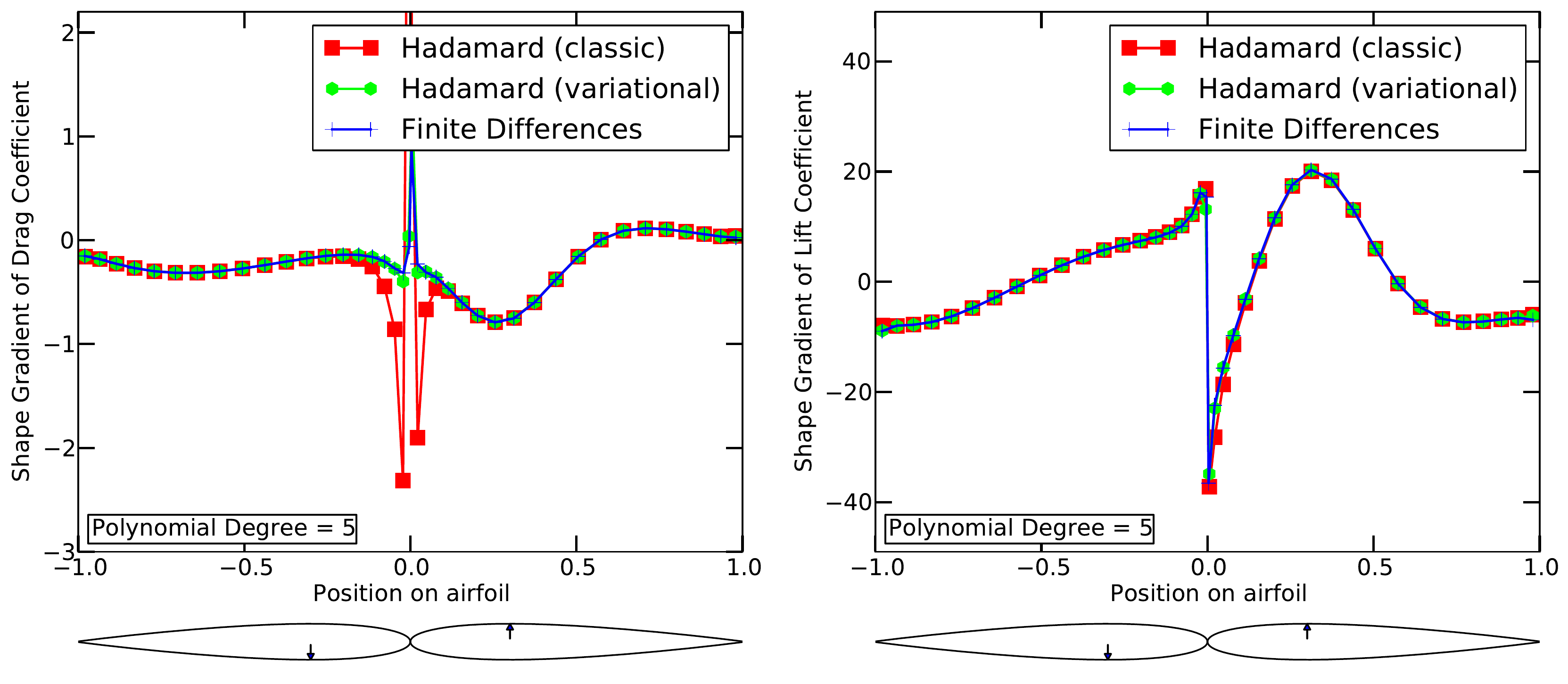}
\caption{Shape Gradient of Drag/Lift Coefficient to polynomial degree $p=5$.}
\label{fi:draglift3}
\end{figure}

Further $p$-refinement and therewith enhancement of the accuracy of the solution follows this trend: In Figure \ref{fi:draglift3} we see an excellent match between the variational approach and finite differences, whereas the pointwise approach still exhibits a fairly substantial gap in the region around the nose, although somewhat diminished when compared to the lower degree case.

The natural next step after deriving the gradients is conducting the actual design optimization. A very popular choice for an optimization strategy is the so called one-shot methodology, where the design update is made simultaneously to the iteration of the primal and adjoint solver~\cite{GS09, GriewankxX2005, SSIG2010A}. Consequently, gradients created from an inexact and not fully converged flow and adjoint solution are used initially, before all three residuals are driven to convergence simultaneously. Thus, the precision of the shape derivative for inaccurate flow solutions is of crucial importance.

Despite the shortcomings in accuracy of the polynomial degree three solution, as shown in Figure \ref{fi:draglift1}, it is well worth pointing out that the quality of the variational Hadamard form would most likely suffice to be used in a one-shot optimization without compunction. This can be justified, because the general manner of the shape derivative, especially in terms of sign, is captured. Contrary to this, the shape derivative of the pointwise approach shows a tremendous deviation over the whole profile, which questions its utility for optimization, which we intend to study as future work.  

\section{Conclusion}
Shape optimization under PDE constraints very often follows the ``function composition'' approach, where the existence of local shape derivatives for each component in the whole chain containing mesh perturbation, PDE variation and the variation of the objective is assumed to exist. For elliptic problems, this approach is known to work very well~\cite{sokolowski}, somewhat contrary to problems within aerodynamic design, where usually great care needs to be taken~\cite{Lozano2012}. The main purpose of the present work was to circumvent this ``step-by-step'' differentiation and directly shape differentiate the weak form of the governing equations, similar to~\cite{KunishVariationalSD}. Here, however, we focus on the challenging task of shape optimization within compressible, viscous fluids, for which we also demonstrate the respective gain in accuracy numerically. Any approach based on the weak form is also much more inline with the actual flow solver. Thus, a variational methodology greatly benefits the alignment of the continuous shape differentiation process with the discrete implementation.

One aspect of the present work was the direct comparison between both approaches, which revealed extra terms arising if the variational form of the state equation is used as the governing model. Both shape derivatives were implemented into the DLR flow solver PADGE, a discontinuous Galerkin flow solver of variable order operating on fourth order curvilinear meshes. We conclude with numerical accuracy studies where we gauge both derivatives against a reference solution created by finite differences. Although the gap between the shape derivative based on the weak and strong form diminishes slightly with increasing order of the polynomial ansatz functions, we always found the weak form derivative to be of considerably higher accuracy. Finally, we found the general quality of this weak form shape derivative to be very promising for a future application within a one-shot optimization framework.

\section*{Acknowledgments}
We wish to thank the BMBF for supporting this work as part of the project DGHPOPT BMBF 05M10PAB. We also wish to thank Dr.\ R. Hartmann, DLR Braunschweig, for his assistance with the flow solver PADGE.  

\bibliographystyle{plain}
\bibliography{SSG_ShapeDerivativesForTheCompressibleNavierStokesEquationsInVariationalForm}

\begin{thebibliography}{10}

\bibitem{Arian2}
E.~Arian and V.~N. Vatsa.
\newblock A preconditioning method for shape optimization governed by the
  {E}uler equations.
\newblock Technical Report 98-14, Institute for Computer Applications in
  Science and Engineering (ICASE), 1998.

\bibitem{castro}
C.~Castro, C.~Lozano, F.~Palacios, and E.~Zuazua.
\newblock Systematic continuous adjoint approach to viscous aerodynamic design
  on unstructured grids.
\newblock {\em AIAA}, 45(9):2125--2139, 2007.

\bibitem{shapesandgeometries}
M.~C. Delfour and J.-P. Zol{\'e}sio.
\newblock {\em Shapes and Geometries}.
\newblock Society for Industrial and Applied Mathematics, Philadelphia, 2001.

\bibitem{Eppler2007A}
K.~Eppler, H.~Harbrecht, and R.~Schneider.
\newblock On convergence in elliptic shape optimization.
\newblock {\em SIAM J. Control Optim.}, 46:61--83, March 2007.

\bibitem{GaWaMoWi07}
N.~Gauger, A.~Walther, C.~Moldenhauer, and M.~Widhalm.
\newblock Automatic differentiation of an entire design chain for aerodynamic
  shape optimization.
\newblock {\em Notes on Numerical Fluid Mechanics and Multidisciplinary
  Design}, 96:454--461, 2007.

\bibitem{giles1997adjoint}
M.~B. Giles and N.~A. Pierce.
\newblock Adjoint equations in {C}{F}{D}: duality, boundary conditions and
  solution behaviour.
\newblock {\em AIAA}, 97-1850, 1997.

\bibitem{GriewankxX2005}
A.~Griewank.
\newblock Projected {H}essians for preconditioning in one-step one-shot design
  optimization.
\newblock {\em Nonconvex Optimization and its Application}, 83:151--172, 2006.

\bibitem{Har08b}
R.~Hartmann.
\newblock Numerical analysis of higher order discontinuous galerkin finite
  element methods.
\newblock In H.~Deconinck, editor, {\em VKI LS 2008-08: CFD - ADIGMA course on
  very high order discretization methods, Oct. 13-17, 2008}. Von Karman
  Institute for Fluid Dynamics, Rhode Saint Gen{\`e}se, Belgium, 2008.

\bibitem{KunishVariationalSD}
K.~Ito, K.~Kunisch, G.~Gunther, and H.~Peichl.
\newblock Variational approach to shape derivatives.
\newblock {\em Control, Optimisation and Calculus of Variations}, 14:517--539,
  2008.

\bibitem{ItoKunischPeichel2006}
K.~Ito, K.~Kunisch, and G.~Peichl.
\newblock Variational approach to shape derivatives of a class of {B}ernoulli
  problems.
\newblock {\em Journal of Mathematical Analysis and Applications},
  314(1):126--149, 2006.

\bibitem{Jameson88}
A.~Jameson.
\newblock Aerodynamic design via control theory.
\newblock {\em Journal of Scientific Computing}, 3(3):233--260, 1988.

\bibitem{jameson1995optimum}
A.~Jameson.
\newblock Optimum aerodynamic design using cfd and control theory.
\newblock volume~97. American Institute of Aeronautics and Astronautics, 1997.

\bibitem{jameson1998optimum}
A.~Jameson, L.~Martinelli, and N.A. Pierce.
\newblock Optimum aerodynamic design using the navier--stokes equations.
\newblock {\em Theoretical and Computational Fluid Dynamics}, 10(1):213--237,
  1998.

\bibitem{Lozano2012}
C.~Lozano.
\newblock Discrete surprises in the computation of sensitivites from boundary
  integrals in the continuous adjoint approach to inviscid aerodynamic shape
  optimization.
\newblock {\em Computers \& Fluids}, 56, 2012.

\bibitem{PironneauxX2001}
B.~Mohammadi and O.~Pironneau.
\newblock {\em Applied Shape Optimization for Fluids}.
\newblock Numerical Mathematics and Scientific Computation. Clarendon Press
  Oxford, 2001.

\bibitem{Zingg:AnaheimJournal}
M.~Nemec and D.~W. Zingg.
\newblock A {N}ewton--{K}rylov algorithm for aerodynamic design using the
  {N}avier--{S}tokes equations.
\newblock {\em AIAA Journal}, 40(6):1146--1154, 2002.

\bibitem{NielsenPark2006}
E.~J. Nielsen and M.~A. Park.
\newblock Using an adjoint approach to eliminate mesh sensitivities in
  computational design.
\newblock {\em AIAA Journal}, 44(5):948--953, 2006.

\bibitem{Papdimitriou2007}
D.~I. Papadimitriou and K.~C. Giannakoglou.
\newblock A continous adjoint method with objective function derivatives based
  on boundary integrals, for inviscid and viscous flows.
\newblock {\em Computers \& Fluids}, 36:325--341, 2007.

\bibitem{Piro73}
O.~Pironneau.
\newblock On optimum profiles in {S}tokes flow.
\newblock {\em Journal of Fluid Mechanics}, 59(1):117--128, 1973.

\bibitem{schmidt_thesis}
S.~Schmidt.
\newblock {\em Efficient Large Scale Aerodynamic Design Based on Shape
  Calculus}.
\newblock PhD thesis, University of Trier, Germany, 2010.

\bibitem{SIGS2011}
S.~Schmidt, C.~Ilic, N.~Gauger, and V.~Schulz.
\newblock Airfoil design for compressible inviscid flow based on shape
  calculus.
\newblock {\em Optimization and Engineering}, 12(3):349--369, 2011.

\bibitem{SSIG2010A}
S.~Schmidt, C.~Ilic, V.~Schulz, and N.~Gauger.
\newblock Three dimensional large scale aerodynamic shape optimization based on
  shape calculus.
\newblock {\em AIAA Journal}, 51(11):2615--2627, 2013.

\bibitem{schmidt:2562}
S.~Schmidt and V.~Schulz.
\newblock Impulse response approximations of discrete shape {H}essians with
  application in {C}{F}{D}.
\newblock {\em SIAM Journal on Control and Optimization}, 48(4):2562--2580,
  2009.

\bibitem{SS10GeneralNS}
S.~Schmidt and V.~Schulz.
\newblock Shape derivatives for general objective functions and the
  incompressible {N}avier-{S}tokes equations.
\newblock {\em Control and Cybernetics}, 39(3):677--713, 2010.

\bibitem{GS09}
V.~Schulz and I.~Gherman.
\newblock One-shot methods for aerodynamic shape optimization.
\newblock In N.~Kroll, D.~Schwamborn, K.~Becker, H.~Rieger, and F.~Thiele,
  editors, {\em MEGADESIGN and MegaOpt --- German Initiatives for Aerodynamic
  Simulation and Optimization in Aircraft Design}, volume 107 of {\em Notes on
  Numerical Fluid Mechanics and Multidisciplinary Design}, pages 207--220.
  Springer, 2009.

\bibitem{Soemarwoto1}
B.~Soemarwoto.
\newblock The variational method for aerodynamic optimization using the
  {N}avier--stokes equations.
\newblock Technical Report 97-71, Institute for Computer Applications in
  Science and Engineering (ICASE), 1997.

\bibitem{sokolowski}
J.~Sokolowski and J.-P. Zol{\'e}sio.
\newblock {\em Introduction to Shape Optimization}.
\newblock Springer-Verlag, Berlin, 1992.

\bibitem{TaAsanSxX1992a}
S.~Ta'asan, G.~Kuruvila, and {M.D.} Salas.
\newblock Aerodynamic design and optimization in one shot.
\newblock In {\em 30th Aerospace Sciences Meeting, Reno, NV, AIAA Paper
  92-0025}, 1992.

\bibitem{Giannakoglou2009B}
A.~S. Zymaris, D.~I. Papadimitriou, K.~C. Giannakoglou, and C.~Othmer.
\newblock Continuous adjoint approach to the {S}palart--{A}llmaras turbulence
  model for incompressible flows.
\newblock {\em Computers \& Fluids}, 38:1528--1538, 2009.

\end{thebibliography}
\end{document}